\documentclass[11pt]{article}

\usepackage{amsmath}
\usepackage{amssymb}
\usepackage[mathcal]{euscript}
\usepackage{graphicx,psfrag}
\usepackage{epsfig}
\usepackage{epic,eepic}

\newcommand{\ee}{{\bf e}}

\newcommand{\bbZ}{{\mathbb Z}}
\newcommand{\bbR}{{\mathbb R}}

\newcommand{\bbC}{{\mathbb C}}
\newcommand{\bbL}{{\mathbb L}}
\newcommand{\bbS}{{\mathbb S}}
\newcommand{\bbP}{{\mathbb P}}

\newcommand{\cC}{{\mathcal C}}

\newcommand{\cQ}{{\mathcal Q}}

\newcommand{\cX}{{\mathcal X}}

\newtheorem{theorem}{Theorem}

\newtheorem{definition}[theorem]{Definition}
\newtheorem{corollary}[theorem]{Corollary}

\begin{document}

\title{Isothermic surfaces in sphere geometries as Moutard nets}
\author{
Alexander I. Bobenko
\thanks{Institut f\"ur Mathematik,
Technische Universit\"at Berlin, Str. des 17. Juni 136, 10623
Berlin, Germany. E--mail: {\tt bobenko@math.tu-berlin.de}} \and
Yuri B. Suris
\thanks{Zentrum Mathematik, Technische Universit\"at M\"unchen,
Boltzmannstr. 3, 85747 Garching bei M\"unchen, Germany. E--mail:
{\tt suris@ma.tum.de}}}

\maketitle{\renewcommand{\thefootnote}{}
\footnote[0]{Research for
this article was supported by the DFG Research Unit 565
`Polyhedral Surfaces'' and the DFG Research Center
\textsc{Matheon} ``Mathematics for key technologies'' in Berlin.}}

{\small {\bf Abstract.}
  We give an elaborated treatment of discrete
isothermic surfaces and their analogs in different geometries
(projective, M\"obius, Laguerre, Lie). We find the core of the
theory to be a novel projective characterization of discrete
isothermic nets as Moutard nets. The latter belong to projective
geometry and are nets with planar faces defined through a
five-point property: a vertex and its four diagonal neighbors span
a three dimensional space. Analytically this property is
equivalent to the existence of representatives in the space of
homogeneous coordinates satisfying the discrete Moutard equation.
Restricting the projective theory to quadrics, we obtain Moutard
nets in sphere geometries.

In particular, Moutard nets in M\"obius geometry are shown to
coincide with discrete isothermic nets. The five-point property in
this particular case says that a vertex and its four diagonal
neighbors lie on a common sphere, which is a novel
characterization of discrete isothermic surfaces. Discrete
Laguerre isothermic surfaces are defined through the corresponding
five-plane property which requires that a plane and its four
diagonal neighbors share a common touching sphere. Equivalently,
Laguerre isothermic surfaces are characterized by having an
isothermic Gauss map. We conclude with Moutard nets in Lie
geometry.}

\newpage
\section{Introduction}
\label{Sect: intro}

This paper is a sequel to our paper ``On organizing principles of
discrete differential geometry. Geometry of spheres'' \cite{BS2},
where the following discretization principles have been
formulated:
\begin{itemize}
 \item {\em Transformation group principle:} smooth geometric objects and
    their discretizations belong to the same geometry, i.e.
    are invariant with respect to the same transformation group.
 \item {\em Consistency principle:} discretizations of smooth
    parametrized geometries can be extended to multidimensional
    consistent nets.
\end{itemize}
Being applied to discretization of curvature line parametrizations
of general surfaces, these principles led to the definition of
{\em principal contact element nets}. These are nets of contact
elements with the property that neighboring contact elements share
a common sphere. In particular, it was shown that the points and
the planes of principal contact element nets build circular and
conical nets, respectively.

In the present paper we turn to isothermic surfaces which is a
special class of surfaces admitting a conformal curvature line
parametrization. This important class of surfaces has been studied
by classics \cite{Da}. In the 1990-s a relation to the theory of
integrable systems has been discovered \cite{CGS, BHPP, BP1}. An
overview of the modern theory of isothermic surfaces can be found
in \cite{HJ}. The theory has been extended for isothermic surfaces
in spaces of arbitrary dimension in \cite{Bu, Sch}. In \cite{BP1}
the theory has been discretized: discrete isothermic surfaces were
defined as special circular nets with factorized cross-ratios of
elementary quadrilaterals. This property is manifestly
M\"obius-invariant. Moreover, it can be consistently imposed on
three-dimensional nets \cite{HHP, BP2}. Thus, discrete isothermic
surfaces is an instance of geometry satisfying both discretization
principles.

In this paper, we give an elaborated treatment of discrete
isothermic surfaces and their analogs in different geometries
(projective, M\"obius, Laguerre, Lie), applying the discretization
principles systematically. We find the core of the theory to be a
novel projective characterization of discrete isothermic nets as
{\em Moutard nets}. The latter belong to projective geometry and
are nets with planar faces $f:\bbZ^2\to\bbR\bbP^N$ defined through
a five-point property: a vertex and its four diagonal neighbors
span a three dimensional space (thus, in comparison with a generic
net with planar faces, the dimension drops by one). Analytically
this property is equivalent to the existence of representatives in
the space of homogeneous coordinates $y:\bbZ^2\to\bbR^{N+1}$
satisfying the {\em Moutard equation}
$$
\tau_1\tau_2 y-y \parallel \tau_2 y-\tau_1 y,
$$
where $\tau_i$ is the shift in the $i$-th coordinate direction.
The consistency of the five-point property follows from the
consistency of the Moutard equation.

Restricting the projective theory to quadrics, we obtain Moutard
nets in sphere geometries. In particular, Moutard nets in M\"obius
geometry are shown to coincide with discrete isothermic nets. The
five-point property in this particular case says that a vertex and
its four diagonal neighbors lie on a common sphere, which is a
novel characterization of discrete isothermic surfaces.

In Laguerre geometry discrete surfaces are maps $\bbZ^2\to\{ {\rm
planes\ in}\ \bbR^3 \}$. The planarity of faces in the Laguerre
quadric is equivalent to the conical property, while the
five-plane property requires that a plane and its four diagonal
neighbors share a common touching sphere. This is a definition of
discrete Laguerre isothermic surfaces. Equivalently, discrete
Laguerre isothermic surfaces are characterized by having an
isothermic Gauss map. The latter class was independently
introduced in \cite{WP}. Smooth Laguerre isothermic surfaces have
been studied in \cite{E, MN1, MN2}.

We conclude with Moutard nets in Lie geometry. The latter are
special Ribaucour sphere congruences $\bbZ^2\to\{{\rm spheres\ in}
\ \bbR^3\}$ with the corresponding five-sphere property. A
particular case is S-isothermic surfaces \cite{BP2, BHS, Ho}.

\section{Discrete Moutard nets}
\label{subsect: moutard}

In considerations of various nets $f:\bbZ^2\to\cX$, we use the
following notational conventions: for some fixed $u\in\bbZ^2$, we
write $f$ for $f(u)$, and further $f_i$ for $\tau_if(u)=f(u+e_i)$,
and $f_{-i}$ for $\tau_i^{-1}f(u)=f(u-e_i)$. Also, we freely use
notations, definitions and results from \cite{BS2}.

\begin{definition}\label{def: dmn}
{\bf (Discrete Moutard net)} A two-dimensional Q-net
$f:\bbZ^2\to\bbR\bbP^{N}$ $(N\ge 4)$ is called a discrete Moutard
net, if for any $u\in\bbZ^2$ the five points $f$ and $f_{\pm 1,\pm
2}$ lie in a three-dimensional subspace $V\subset\bbR\bbP^{N}$,
not containing some (and then any) of the four points $f_{\pm 1}$,
$f_{\pm 2}$.
\end{definition}
Thus, the defining condition of a discrete Moutard net deals with
four elementary planar quadrilaterals adjacent to one vertex. As a
consequence of this definition, all nine vertices of such four
quadrilaterals of a discrete Moutard net lie in a four-dimensional
subspace of $\bbR\bbP^N$.

\begin{theorem}\label{Thm: charact 2d T-net}
{\bf (Discrete Moutard equation)} A discrete Moutard net
$\,f:\bbZ^2\to\bbR\bbP^{N}$ possesses a lift to the homogeneous
coordinates space $y:\bbZ^2\to\bbR^{N+1}$ satisfying the discrete
Moutard equation:
\begin{equation}\label{eq:dMou 2d minus}
y_{12}-y=a_{12} (y_2-y_1),
\end{equation}
with some $a_{12}:\bbZ^2\to\bbR$ (it is natural to assign the real
numbers $a_{12}$ to the elementary squares of $\,\bbZ^2$).
\end{theorem}
{\bf Proof.} We start with the observation that for {\em any}
Q-net $f$ in $\bbR\bbP^{N}$ it is always possible (and almost
trivial) to find homogeneous coordinates for the four vertices of
{\em one} elementary quadrilateral satisfying the discrete Moutard
equation on that quadrilateral. Moreover, one can do this for an
arbitrary choice of homogeneous coordinates for any two
neighboring vertices of the quadrilateral. Indeed, consider any
homogeneous coordinates
$\tilde{f},\tilde{f}_1,\tilde{f}_2,\tilde{f}_{12}\in\bbR^{N+1}$
for the vertices of a planar quadrilateral, connected by a linear
relation
\[
\tilde{f}_{12}=\tilde{c}_{21}\tilde{f}_1+
\tilde{c}_{12}\tilde{f}_2+\rho_{12}\tilde{f}.
\]
Keep the representatives $y=\tilde{f}$, $y_1=\tilde{f}_1$, say,
and set $\tilde{f}_{12}=\rho_{12}y_{12}$ and $\tilde{f}_2=ay_2$
with $a=-\tilde{c}_{21}/\tilde{c}_{12}$, then $y$ satisfies the
discrete Moutard equation (\ref{eq:dMou 2d minus}) within one
elementary quadrilateral.

Now, for Q-nets with a special property formulated in Definition
\ref{def: dmn}, this construction can be extended to the whole
net. Start with arbitrary representatives $y$, $y_1$, and proceed
clockwise around the vertex $y$. We find consecutively: the
representatives $y_{-2}$, $y_{1,-2}$ which assure the Moutard
equation on the quadrilateral $(y,y_1,y_{1,-2},y_{-2})$, then the
representatives $y_{-1}$, $y_{-1,-2}$ which assure the Moutard
equation on the quadrilateral $(y,y_{-1},y_{-1,-2},y_{-2})$, and
then the representatives $y_{2}$, $y_{-1,2}$ which assure the
Moutard equation on the quadrilateral $(y,y_{-1},y_{-1,2},y_{2})$,
see Fig.~\ref{Fig: towards Moutard}.
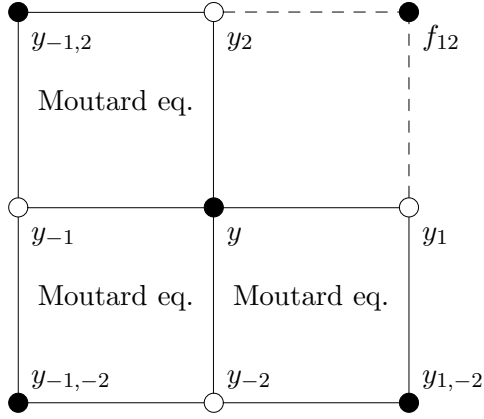
\begin{figure}[htbp]
\begin{center}
\setlength{\unitlength}{0.045em}
\begin{picture}(300,320)(0,0)
 \put(0,0){\circle*{15}} \put(150,0){\circle{15}}
 \put(300,0){\circle*{15}}
 \put(0,150){\circle{15}} \put(150,150){\circle*{15}}
 \put(300,150){\circle{15}}
 \put(0,300){\circle*{15}} \put(150,300){\circle{15}}
 \put(300,300){\circle*{15}}
 \path(0,0)(142.5,0) \path(157.5,0)(300,0)
 \path(0,0)(0,142.5) \path(0,157.5)(0,300)
 \path(150,7.5)(150,292.5)
 \path(7.5,150)(292.5,150)
 \path(0,300)(142.5,300)    \path(300,0)(300,142.5)
 \dashline[+30]{10}(157.5,300)(300,300)
 \dashline[+30]{10}(300,157.5)(300,300)
 \put(160,125){$y$}
 \put(310,125){$y_1$}
 \put(10,125){$y_{-1}$}
 \put(160,15){$y_{-2}$}
 \put(310,15){$y_{1,-2}$}
 \put(10,15){$y_{-1,-2}$}
 \put(160,275){$y_2$}
 \put(10,275){$y_{-1,2}$}
 \put(310,275){$f_{12}$}
 \put(15,75){Moutard eq.}
 \put(15,225){Moutard eq.}
 \put(165,75){Moutard eq.}
\end{picture}
\caption{Constructing a Moutard representative for a projective
Q-net with a three-dimensional black cross}\label{Fig: towards
Moutard}
\end{center}
\end{figure}
In the remaining quadrilateral $(y,y_1,y_{12},y_2)$, the
representatives $y$, $y_1$, $y_2$ are already fixed on the
previous steps of the construction, so that we can dispose of the
representative $y_{12}$ of $f_{12}$ only. Observe that the point
with the representative $y_1-y_2$ belongs to the the plane
$\Pi\subset\bbR\bbP^N$ of the quadrilateral $(f,f_1,f_{12},f_2)$
(obviously) and to the three-dimensional space
$V\subset\bbR\bbP^N$ through the points $f$, $f_{1,-2}$,
$f_{-1,-2}$, $f_{-1,2}$, because of the equation
\begin{eqnarray*}
y_1-y_2 & = & (y_1-y_{-2})+(y_{-2}-y_{-1})+(y_{-1}-y_2)\\
 & = & \alpha(y_{1,-2}-y)+\beta(y_{-1,-2}-y)+\gamma(y_{-1,2}-y).
 \end{eqnarray*}
By the hypothesis of the theorem, the point $f_{1,2}$ lies in the
latter space $V$. Therefore, the whole line through $f$ and
$f_{1,2}$ lies in the intersection $\Pi\cap V$. Since $N\ge 4$, we
conclude that in general position $\Pi\cap V$ {\em is} the line
through $f$ and $f_{12}$. Thus, the point with the representative
$y_1-y_2$ belongs to this line, therefore $y_1-y_2$ is a linear
combination of $y$ and $y_{12}$. By a suitable choice of
representative $y_{12}$ of $f_{12}$, we can make $y_1-y_2$
proportional to $y_{12}-y$. Thus, the construction of
representatives satisfying the Moutard equation closes up around
any vertex. This allows to extend the construction to the whole
lattice $\bbZ^2$. $\Box$
\smallskip

Definition \ref{def: dmn} is non-applicable in the case when some,
and then all of the points $f_{\pm 1}$, $f_{\pm 2}$ lie in the
three-dimensional space $V$ through $f$, $f_{\pm 1,\pm 2}$, in
particular, it cannot be used to define discrete Moutard nets in
$\bbR\bbP^3$. We will show that Theorem \ref{Thm: charact 2d
T-net} remains valid if one defines discrete Moutard nets in
$\bbR\bbP^3$ as follows.
\begin{definition}\label{def: dmn 3d}
{\bf (Discrete Moutard net in $\bbR\bbP^3$)} A two-dimensional
Q-net $f:\bbZ^2\to\bbR\bbP^3$ is called a discrete Moutard net, if
for any $u\in\bbZ^2$ the following condition is satisfied: the
three planes
\[
\Pi^{(\rm up)}=(f,f_{12},f_{-1,2}),\quad \Pi^{(\rm
down)}=(f,f_{1,-2},f_{-1,-2}),\quad \Pi^{(1)}=(f,f_1,f_{-1})
\]
have a common line $\ell^{(1)}$.
\end{definition}

{\bf Remark 1.} It is not difficult to see that in the context of
Definition \ref{def: dmn} with $N\ge 4$ the requirement of
Definition \ref{def: dmn 3d} is automatically satisfied. Indeed,
in this case all nine points $f$, $f_{\pm 1}$, $f_{\pm 2}$ and
$f_{\pm 1,\pm 2}$ lie in a four-dimensional subspace of
$\bbR\bbP^N$. In this subspace one can consider, along with the
three-dimensional subspace $V$, the three-dimensional subspaces
$V^{\rm (up)}$ containing the two quadrilaterals
$(f,f_1,f_{12},f_2)$, $(f,f_{-1},f_{-1,2},f_2)$, and $V^{\rm
(down)}$ containing the quadrilaterals $(f,f_1,f_{1,-2},f_{-2})$,
$(f,f_{-1},f_{-1,-2},f_{-2})$. Obviously, one has:
\[
 \Pi^{\rm (up)}=V^{\rm (up)}\cap V,\quad
 \Pi^{\rm (down)}=V^{\rm (down)}\cap V,\quad
 \Pi^{(1)}=V^{\rm (up)}\cap V^{\rm (down)}.
\]
Generically, three three-dimensional subspaces $V$, $V^{\rm (up)}$
and $V^{\rm (down)}$ of a four-dimensional space intersect along a
line $\ell^{(1)}$.

\begin{figure}[htbp]
 \psfrag{f_1}[Bl][bl][0.9]{$f_1$}
 \psfrag{f_2}[Bl][bl][0.9]{$f_2$}
 \psfrag{f_{-1}}[Bl][bl][0.9]{$f_{-1}$}
 \psfrag{f_{-2}}[Bl][bl][0.9]{$f_{-2}$}
 \psfrag{f_{12}}[Bl][bl][0.9]{$f_{12}$}
 \psfrag{f_{1,-2}}[Bl][bl][0.9]{$f_{1,-2}$}
 \psfrag{f_{-1,2}}[Bl][bl][0.9]{$f_{-1,2}$}
 \psfrag{f_{-1,-2}}[Bl][bl][0.9]{$f_{-1,-2}$}
 \psfrag{l^1}[Bl][bl][0.9]{$\ell^{(1)}$}
 \psfrag{l^2}[Bl][bl][0.9]{$\ell^{(2)}$}
 \center{\includegraphics[width=120mm]{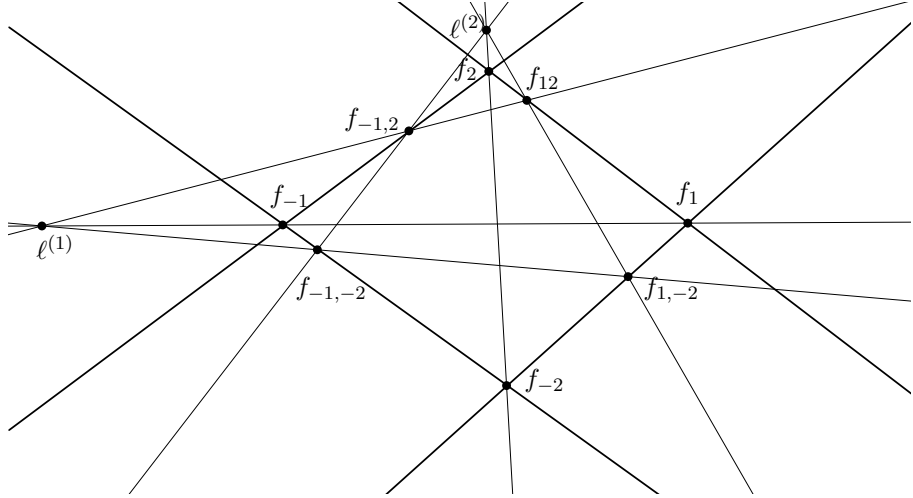}}
 \caption{Desargues theorem}
 \label{Fig: Desargues}
\end{figure}

{\bf Remark 2.} There is an asymmetry between the coordinate
directions 1 and 2 in Definition \ref{def: dmn 3d}. However, this
asymmetry is apparent: the condition in Definition \ref{def: dmn
3d} is equivalent to the requirement that the three planes
\[
\Pi^{(\rm left)}=(f,f_{-1,2},f_{-1,-2}),\quad \Pi^{(\rm
right)}=(f,f_{1,2},f_{1,-2}),\quad \Pi^{(2)}=(f,f_{2},f_{-2})
\]
have a common line $\ell^{(2)}$. One way to see this is to
consider a central projection of the whole picture from the point
$f$ to some plane not containing $f$. In this projection, the
planarity of elementary quadrilaterals $(f,f_i,f_{ij},f_j)$ turns
into collinearity of the triples of points $f_i$, $f_j$ and
$f_{ij}$. The traces of the planes $\Pi^{(\rm up)}$, $\Pi^{(\rm
down)}$ and $\Pi^{(1)}$ on the projection plane are the lines
$(f_{12},f_{-1,2})$, $(f_{1,-2},f_{-1,-2})$, and $(f_1,f_{-1})$,
respectively, and the requirement of Definition \ref{def: dmn 3d}
turns into the requirement for these three lines to meet in a
point. Similarly, the traces of the planes $\Pi^{(\rm left)}$,
$\Pi^{(\rm right)}$ and $\Pi^{(2)}$ on the projection plane are
the lines $(f_{-1,2},f_{-1,-2})$, $(f_{1,2},f_{1,-2})$, and
$(f_2,f_{-2})$, respectively. The requirement for the latter three
lines to meet in a point is equivalent to the  previous one  --
this is the statement of the famous Desargues theorem, see Fig.
\ref{Fig: Desargues}.

Another way to demonstrate the actual symmetry between the
coordinate directions 1 and 2 in Definition \ref{def: dmn 3d} is
to show that Theorem \ref{Thm: charact 2d T-net} still holds in
$\bbR\bbP^3$. Indeed, the discrete Moutard equation (\ref{eq:dMou
2d minus}) is manifestly symmetric with respect to the flip
$1\leftrightarrow 2$.

{\bf Proof of Theorem \ref{Thm: charact 2d T-net} for $N=3$.} We
start the proof exactly as in the general case $N\ge 4$. The only
thing to be changed is the demonstration of the fact that the
point with the representative $y_2-y_1$ lies on the line through
$f$ and $f_{12}$. To do this in the present situation, we first
observe that, due to
\[
y_1-y_{-1}=(y_1-y_{-2})+(y_{-2}-y_{-1})=
\alpha(y_{1,-2}-y)+\beta(y_{-1,-2}-y),
\]
the point with the representative $y_1-y_{-1}$ lies in the plane
$\Pi^{\rm(down)}$. Obviously, it lies also in $\Pi^{(1)}$.
Therefore, it lies on the line $\ell^{(1)}$. As a consequence of
the property of Definition \ref{def: dmn 3d}, it belongs also to
the plane $\Pi^{\rm(up)}$. Now from
\[
y_2-y_1=(y_2-y_{-1})-(y_1-y_{-1})=\gamma(y_{-1,2}-y)-(y_1-y_{-1})
\]
we find that the point with the representative $y_2-y_1$ belongs
to $\Pi^{\rm(up)}$, as well. Since the point with the
representative $y_2-y_1$ also belongs (obviously) to the plane of
the quadrilateral $(f,f_1,f_{12},f_2)$, we conclude that it lies
in the intersection of the latter plane with
$\Pi^{\rm(up)}=(f,f_{12},f_{-1,2})$, which is, in the generic
case, the line through $f$ and $f_{12}$. $\Box$
\medskip

Definitions \ref{def: dmn}, \ref{def: dmn 3d} are essentially
dealing with two-dimensional Q-nets. However, the characterization
of discrete Moutard nets given in Theorem \ref{Thm: charact 2d
T-net} opens a way to define multi-dimensional Moutard nets, and,
in particular, to define transformations of Moutard nets with
remarkable permutability properties. Namely, it turns out that eq.
(\ref{eq:dMou 2d minus}) can be posed on multidimensional
lattices.

\begin{definition}\label{dfn:tn}
{\bf (T-net)} A map $y:\bbZ^m\to\bbR^N$ is called an
$m$-dimensional T-net (trapezoidal net), if for any $u\in\bbZ^m$
and for any pair of indices $i\neq j$ there holds the discrete
Moutard equation
\begin{equation}\label{eq:dMou}
 y_{ij}-y=a_{ij} (y_j-y_i),
\end{equation}
with some $a_{ij}:\bbZ^m\to\bbR$, in other words, if all
elementary quadrilaterals $(y,y_i,y_{ij},y_j)$ are planar and have
parallel diagonals.
\end{definition}
Of course, coefficients $a_{ij}$ have to be skew-symmetric,
$a_{ij}=-a_{ji}$. We show that three-dimensional T-nets are
described by a well-defined three-dimensional system.
\begin{theorem} \label{Th: T-cube}
{\bf (Elementary hexahedron of a T-net)} Given seven points $y$,
$y_i$, and $y_{ij}$ $\,(1\leq i\neq j\leq 3)$ in $\,\bbR^N$, such
that eq. (\ref{eq:dMou}) is satisfied on the three quadrilaterals
$(y,y_i,y_{ij},y_j)$ adjacent to the vertex $y$, there exists a
unique point $y_{123}$ such that eq. (\ref{eq:dMou}) is satisfied
on the three quadrilaterals $(y_i,y_{ij},y_{123},y_{ik})$ adjacent
to the vertex $y_{123}$.
\end{theorem}
{\bf Proof.} Three equations (\ref{eq:dMou}) for the faces of an
elementary cube of $\bbZ^3$ adjacent to $y_{123}$, give:
\[
\tau_iy_{jk}=\big(1+(\tau_ia_{jk})(a_{ij}+a_{ki})\big)y_i
-(\tau_ia_{jk})a_{ij}y_j-(\tau_ia_{jk})a_{ki}y_k.
\]
They lead to consistent results for $y_{123}$ for arbitrary
initial data, if and only if the following conditions are
satisfied:
\begin{eqnarray*}
1+(\tau_1a_{23})(a_{12}+a_{31}) & = & -(\tau_2a_{31})a_{12}\;=\;
-(\tau_3a_{12})a_{31},\\
1+(\tau_2a_{31})(a_{23}+a_{12}) & = & -(\tau_3a_{12})a_{23}\;=\;
-(\tau_1a_{23})a_{12},\\
1+(\tau_3a_{12})(a_{23}+a_{31}) & = & -(\tau_1a_{23})a_{31}\;=\;
-(\tau_2a_{31})a_{23}.
\end{eqnarray*}
These conditions constitute a system of 6 (linear) equations for 3
unknown variables $\tau_i a_{jk}$ in terms of the known ones
$a_{jk}$. A direct computation shows that this system is not
overdetermined but admits a unique solution:
\begin{equation}\label{eq:star-triang 3}
\frac{\tau_1a_{23}}{a_{23}}=\frac{\tau_2a_{31}}{a_{31}}=
\frac{\tau_3a_{12}}{a_{12}}=
-\frac{1}{a_{12}a_{23}+a_{23}a_{31}+a_{31}a_{12}}\,.
\end{equation}
With $\tau_ia_{jk}$ so defined, eqs. (\ref{eq:dMou}) are fulfilled
on all three quadrilaterals adjacent to $y_{123}$. $\Box$

Eqs. (\ref{eq:star-triang 3}) represent a well-defined birational
map $\{a_{jk}\}\mapsto\{\tau_ia_{jk}\}$, which can be considered
as the fundamental 3D system related to T-nets. It is sometimes
called the {\em ``star-triangle map''}.

Theorem \ref{Th: T-cube} means that the defining condition of
T-nets (parallel diagonals of elementary planar quadrilaterals)
yields a discrete 3D system with fields on vertices taking values
in an affine space $\bbR^N$. This system can be considered as an
admissible reduction of the 3D system describing Q-nets in
$\bbR^N$. Indeed, if one has an elementary hexahedron of an affine
Q-net $y:\bbZ^3\to\bbR^N$ such that its elementary quadrilaterals
$(y,y_i,y_{ij},y_j)$ have parallel diagonals, then the elementary
quadrilaterals $(y_i,y_{ij},y_{123},y_{ik})$ have this property,
as well. To see this, observe that the point $y_{123}$ from
Theorem \ref{Th: T-cube} satisfies the planarity condition, and
therefore it has to coincide with the unique point defined by
planarity of the quadrilaterals $(y_i,y_{ij},y_{123},y_{ik})$.
\smallskip

The 4D consistency of T-nets is a consequence of the analogous
property of Q-nets, since T-constraint propagates in the
construction of a Q-net from its coordinate surfaces. On the level
of formulas we have for T-nets with $m\ge 4$ the system
(\ref{eq:dMou}), while the map $\{a_{jk}\}\mapsto\{\tau_ia_{jk}\}$
is given by
\begin{equation}\label{eq:star-triang}
\frac{\tau_ia_{jk}}{a_{jk}}=
-\frac{1}{a_{ij}a_{jk}+a_{jk}a_{ki}+a_{ki}a_{ij}}\,.
\end{equation}
All indices $i,j,k$ vary now between 1 and $m$, and for any triple
of pairwise different indices $(i,j,k)$, equations involving these
indices solely, form a closed subset.
\medskip

The multidimensional consistency of T-nets yields in a usual
fashion Darboux transformations with permutability properties
(which in the present context should be called discrete Moutard
transformations). We refer to \cite{BS1} for the background on the
relation of multidimensional consistency to Darboux
transformations, and give here only the formulas for the discrete
Moutard transformation of eq. (\ref{eq:dMou}) into
\begin{equation}\label{eq:dMou transf}
 y_{ij}^+-y^+=a_{ij}^+ (y_j^+-y_i^+).
\end{equation}
These formulas read:
\begin{equation}\label{eq:dmp T}
y^+_i-y=b_i(y^+-y_i),
\end{equation}
where the quantities $b_i$ and the transformed coefficients
$a_{ij}^+$ are defined by equations
\begin{equation}\label{eq:dMou T}
\frac{\tau_i b_j}{b_j}=
\frac{a_{ij}^+}{a_{ij}}=\frac{1}{(b_i-b_j)a_{ij}+b_ib_j}.
\end{equation}
It is not difficult to recognize in eqs. (\ref{eq:dMou transf}),
(\ref{eq:dmp T}) the same Moutard equations (\ref{eq:dMou}) on the
$(m+1)$-dimensional lattice, with the superscript ``$+$'' used to
denote the shift $\tau_{m+1}$. Similarly, eqs. (\ref{eq:dMou T})
are nothing but the star-triangle formulas (\ref{eq:star-triang})
with $b_i=a_{i,m+1}$.
\medskip

{\bf Remark.} We learned about the projective five-point
characterization of two-dimensional Moutard nets from
conversations with A.~Doliwa. It should be noted that three- and
higher-dimensional Moutard nets admit a different projective
characterization (planarity of tetrahedra formed by odd or by even
vertices of any elementary cube), see \cite{Do3}. Also the paper
\cite{Do2} by A.~Doliwa deals with a closely related notion of
discrete Koenigs nets.

\subsection{Discrete Moutard nets in quadrics}
\label{Sect: dMou in quadric}

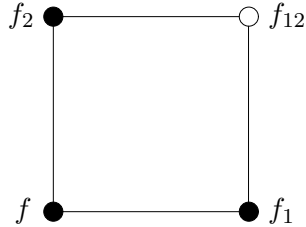
\begin{figure}[htbp]
\begin{center}
\setlength{\unitlength}{0.045em}
\begin{picture}(160,160)(0,0)
 \put(0,0){\circle*{15}}    \put(150,0){\circle*{15}}
 \put(0,150){\circle*{15}}  \put(150,150){\circle{15}}
 \path(0,0)(150,0)       \path(0,0)(0,150)
 \path(150,0)(150,142.5)   \path(0,150)(142.5,150)
 \put(-30,-5){$f$}
 \put(-35,145){$f_2$}  \put(165,-5){$f_1$}
 \put(165,145){$f_{12}$}
\end{picture}
\caption{2D system on an elementary quadrilateral}
\label{Fig:square eq}
\end{center}
\end{figure}

We have seen that discrete Moutard nets (or, more precisely, their
T-net representatives) constitute an admissible reduction of
Q-nets. The restriction to a quadric constitutes another
admissible reduction \cite{Do1}. Imposing two admissible
reductions simultaneously, one comes to T-nets in quadrics. Let
$\bbR^N$ be equipped with a non-degenerate symmetric bilinear form
$\langle\cdot,\cdot\rangle$ (which does not need to be
positive-definite), and let
\[
\cQ=\{y\in\bbR^N: \langle y,y\rangle=\kappa_0\}
\]
be a quadric in $\bbR^N$. We study T-nets $y:\bbZ^m\to\cQ$. This
leads to a {\em discrete 2D system}, since constructing elementary
quadrilaterals of T-nets in $\cQ$ corresponding to elementary
squares of the lattice $\bbZ^m$ admits a well-posed initial value
problem: given three points $y,y_1,y_2\in \cQ$, one finds a unique
fourth point $y_{12}\in\cQ,\,$ $y_{12}\neq y,\,$ satisfying the
discrete Moutard equation
\[
y_{12}-y=a_{12}(y_2-y_1).
\]
Indeed, the condition
\[
\langle y_{12},y_{12}\rangle=\langle
y+a_{12}(y_2-y_1),y+a_{12}(y_2-y_1)\rangle=\kappa_0
\]
leads to a quadratic equation for $a_{12}$, which has one trivial
solution $a_{12}=0\;\Leftrightarrow\; y_{12}=y$, and one
non-trivial:
\[
a_{12}=\frac{\langle y,y_1-y_2\rangle}
 {\kappa_0-\langle y_1,y_2\rangle}\,.
\]
This elementary construction step, i.e., finding the fourth vertex
of an elementary quadrilateral out of the known three vertices, is
symbolically represented on Fig.~\ref{Fig:square eq}.

Turning to an elementary cube of dimension $m\ge 3$, we see that
one can prescribe all points $y$ and $y_i$ for all $1\leq i\leq
m$. Indeed, these data are independent, and one can construct all
other vertices of an elementary cube from these data, {\em
provided one does not encounter contradictions}. To see the
possible source of contradictions, consider in detail the case of
$m=3$. From $y$ and $y_i$ ($1\le i\le 3$) one determines all
$y_{ij}$ by
\begin{equation}\label{eq: dMou in quadric}
y_{ij}-y=a_{ij}(y_j-y_i),\qquad a_{ij}=\frac{\langle
y,y_i-y_j\rangle}{\kappa_0-\langle y_i,y_j\rangle}\,
\end{equation}
After that one has, in principle, three different ways to
determine $y_{123}$, from three squares adjacent to this point;
see Fig.~\ref{cube}. These three values for $y_{123}$ have to
coincide, independently of initial conditions.
\begin{definition} {\bf (3D consistency)}
A 2D system is called 3D consistent, if it can be imposed on all
two-dimensional faces of an elementary cube of $\bbZ^3$.
\end{definition}
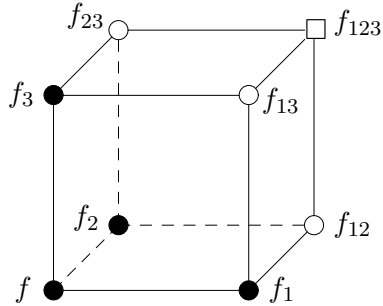
\begin{figure}[htbp]
\begin{center}
\setlength{\unitlength}{0.045em}
\begin{picture}(200,210)(0,0)
 \put(0,0){\circle*{15}}    \put(150,0){\circle*{15}}
 \put(0,150){\circle*{15}}  \put(150,150){\circle{15}}
 \put(50,50){\circle*{15}} \put(50,200){\circle{15}}
 \put(200,50){\circle{15}}
 \put(193,193){$\Box$}
 \path(0,0)(150,0)       \path(0,0)(0,150)
 \path(150,0)(150,142.5)   \path(0,150)(142.5,150)
 \path(0,150)(44.7,194.7)    \path(155.3,155.3)(193,193)
 \path(57.5,200)(193,200)
 \path(200,193)(200,57.5) \path(150,0)(194.7,44.7)
 \dashline[+30]{10}(0,0)(50,50)
 \dashline[+30]{10}(50,50)(50,192.5)
 \dashline[+30]{10}(50,50)(192.5,50)
 \put(-30,-5){$f$}
 \put(-35,145){$f_3$} \put(215,45){$f_{12}$}
 \put(165,-5){$f_1$} \put(160,140){$f_{13}$}
 \put(15,50){$f_2$}  \put(10,205){$f_{23}$}
 \put(215,200){$f_{123}$}
\end{picture}
\caption{3D consistency of 2D systems}\label{cube}
\end{center}
\end{figure}
There holds a quite general theorem, analogous to Theorem 5 of
\cite{BS2}:
\begin{theorem}\label{Thm: 3D yields MD}
{\bf (3D consistency yields consistency in all higher dimensions)}
Any 3D consistent discrete 2D system is also $m$-dimensionally
consistent for all $m>3$.
\end{theorem}
{\bf Proof} goes by induction in $m$ and is analogous to the proof
of Theorem 5 from \cite{BS2}. $\Box$
\smallskip

\begin{theorem}\label{Th: M-lat in quadric}
{\bf (T-nets in quadrics are 3D consistent)} The 2D system
(\ref{eq: dMou in quadric}) governing T-nets in $\cQ$ is 3D
consistent.
\end{theorem}
{\bf Proof.} This can be checked by a tiresome computation, which
can be however avoided by the following conceptual argument.
T-nets in $\cQ$ are a result of imposing two admissible reductions
on Q-nets in $\bbR^N$, namely the T-reduction and the restriction
to a quadric $\cQ$. This reduces the effective dimension of the
system by 1 (allows to determine the fourth vertex of an
elementary quadrilateral from the three known ones), and transfers
the original 3D equation into the 3D consistency of the reduced 2D
equation. Indeed, after finding $y_{12}$, $y_{23}$ and $y_{13}$,
one can construct $y_{123}$ according to the planarity condition
(as intersection of three planes). Then both the T-condition and
the $\cQ$-condition are fulfilled for all three quadrilaterals
adjacent to $y_{123}$, according to Theorem \ref{Th: T-cube} and
the result of \cite{Do1}. Therefore, these quadrilaterals satisfy
our 2D system. $\Box$
\medskip

We mention also an important property of T-nets in quadrics used
in the sequel: the functions
\begin{equation}\label{eq: alpha dMou quadr}
\alpha_i=\langle y,y_i\rangle,
\end{equation}
defined on edges of $\,\bbZ^m$ parallel to the $i$-th coordinate
axes, satisfy
\begin{equation}\label{eq: labelling}
\tau_i\alpha_j=\alpha_j,\qquad i\neq j,
\end{equation}
i.e, any two opposite edges of any elementary square carry the
same value of the corresponding $\alpha_i$. Indeed, equations
\[
\langle y_{ij},y_j\rangle=\langle y_i,y\rangle, \quad \langle
y_{ij},y_i\rangle=\langle y_j,y\rangle.
\]
follow from (\ref{eq: dMou in quadric}) by a direct computation.

\section{Isothermic surfaces in M\"obius geometry}
\label{subsect: Moeb isothermic}

\begin{definition}\label{def: Moeb isothermic}
{\bf (Discrete isothermic surface)} A two-dimensional circular net
$f:\bbZ^2\to\bbR^N$ is called a discrete isothermic surface, if
the corresponding $\hat{f}:\bbZ^2\to\bbL^{N+1,1}$ is a discrete
Moutard net.
\end{definition}

From Definitions \ref{def: dmn}, \ref{def: dmn 3d} there follows a
geometric characterization of discrete isothermic nets:
\begin{theorem}\label{thm: isothermic 5-point sphere}
{\bf (Central spheres for discrete isothermic nets)}

(a) A circular net $f:\bbZ^2\to\bbR^N$ not lying on a
two-dimensional sphere is a discrete isothermic net, if and only
if for any $u\in\bbZ^2$ the five points $f$ and $f_{\pm 1,\pm 2}$
lie on a two-dimensional sphere not containing some (and then any)
of the four points $f_{\pm 1}$, $f_{\pm 2}$.

(b) A circular net $f:\bbZ^2\to\bbS^2\subset\bbR^N$ is a discrete
isothermic net, if and only if for any $u\in\bbZ^2$ the three
circles through $f$,
\begin{eqnarray*}
& C^{\rm(up)}={\rm circle}(f,f_{12},f_{-1,2}),\quad
C^{\rm(down)}={\rm circle}(f,f_{1,-2},f_{-1,-2}), &
\\
& C^{\rm(1)}={\rm circle}(f,f_{1},f_{-1}), &
\end{eqnarray*}
have one additional point in common, which is also equivalent for
the three circles through $f$,
\begin{eqnarray*}
& C^{\rm(left)}={\rm circle}(f,f_{-1,2},f_{-1,-2}),\quad
C^{\rm(right)}={\rm circle}(f,f_{1,2},f_{1,-2}), &
\\
& C^{\rm(2)}={\rm circle}(f,f_{2},f_{-2}), &
\end{eqnarray*}
to have one additional point in common.
\end{theorem}

\begin{figure}[htbp]
\begin{center}
\rotatebox{-90}{\includegraphics[width=0.5\textwidth]{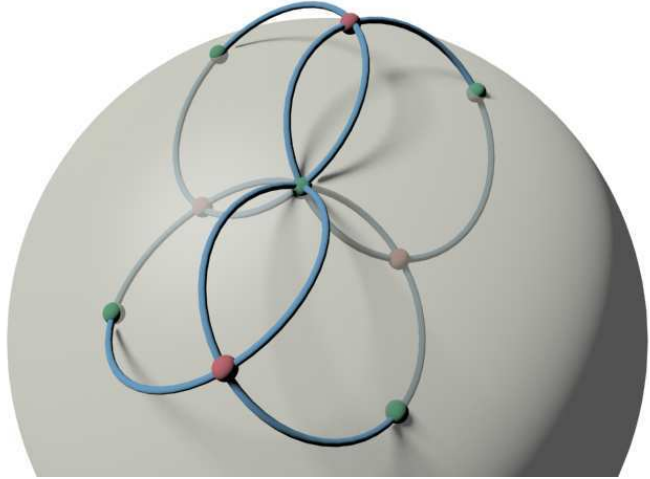}}
\end{center}
\caption{Four circles of a generic discrete isothermic surface,
with a central sphere.} \label{fig: isothermic}
\end{figure}

\begin{figure}[htbp]
\psfrag{f}[Bl][bl][0.9]{$f$}
 \psfrag{f_1}[Bl][bl][0.9]{$f_1$}
 \psfrag{f_2}[Bl][bl][0.9]{$f_2$}
 \psfrag{f_{-1}}[Bl][bl][0.9]{$f_{-1}$}
 \psfrag{f_{-2}}[Bl][bl][0.9]{$f_{-2}$}
 \psfrag{f_{1,2}}[Bl][bl][0.9]{$f_{12}$}
 \psfrag{f_{1,-2}}[Bl][bl][0.9]{$f_{1,-2}$}
 \psfrag{f_{-1,2}}[Bl][bl][0.9]{$f_{-1,2}$}
 \psfrag{f_{-1,-2}}[Bl][bl][0.9]{$f_{-1,-2}$}
 \psfrag{C^up}[Bl][bl][0.9]{$C^{\rm(up)}$}
 \psfrag{C^down}[Bl][bl][0.9]{$C^{\rm(down)}$}
 \psfrag{C^1}[Bl][bl][0.9]{$C^{\rm(1)}$}
 \center{\includegraphics[width=120mm]{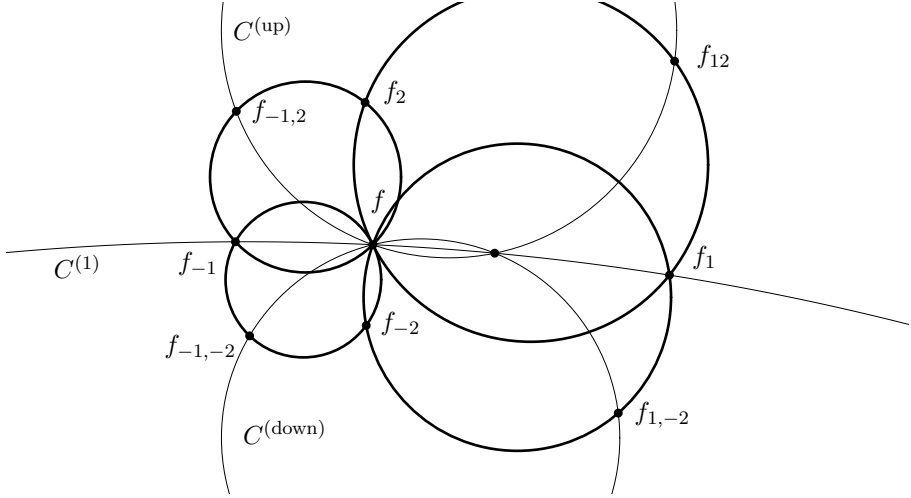}}
 \caption{Four circles of a planar discrete isothermic net.}
 \label{fig: isothermic on a sphere}
\end{figure}

The cases a), b) of Theorem \ref{thm: isothermic 5-point sphere}
are illustrated on Figs. \ref{fig: isothermic}, \ref{fig:
isothermic on a sphere}, respectively.

Another characterization of discrete isothermic surfaces can be
given in terms of the cross-ratios. Recall that for any four
concircular points $f, f_1, f_2, f_{12}\in\bbR^N$ their
(real-valued) cross-ratio can be defined as
\begin{equation}\label{eq: cross-ratio}
q(f,f_1,f_{12},f_2)=(f_1-f)(f_{12}-f_1)^{-1}(f_{12}-f_2)(f_2-f)^{-1}.
\end{equation}
Here multiplication is interpreted as the Clifford multiplication
in the Clifford algebra $\cC\ell(\bbR^N)$. Recall that for
$x,y\in\bbR^N$ the Clifford product satisfies $xy+yx=-2\langle
x,y\rangle$, and that the inverse element of $x\in\bbR^N$ in the
Clifford algebra is given by $x^{-1}=-x/|x|^2$. Alternatively, one
can identify the plane of the quadrilateral $(f,f_1,f_{12},f_2)$
with the complex plane $\bbC$, and then interpret multiplication
in eq. (\ref{eq: cross-ratio}) as the complex multiplication. An
important property of the cross-ratio is its invariance under
M\"obius transformations.
\begin{theorem}\label{thm: isothermic cross-ratios}
{\bf (Four cross-ratios of a discrete isothermic net)} A circular
net $f:\bbZ^2\to\bbR^N$ is a discrete isothermic net, if and only
if the cross-ratios $q=q(f,f_1,f_{12},f_2)$ of its elementary
quadrilaterals satisfy the following condition:
\begin{equation}\label{eq:i/j prod}
    q\cdot q_{-1,-2}=q_{-1}\cdot q_{-2}.
\end{equation}
Here, like in Sect. \ref{subsect: moutard}, the negative indices
$-i$ are used to denote the backward shifts $\tau_i^{-1}$, so
that, e.g., $q_{-1}=q(f_{-1},f,f_2,f_{-1,2})$, see Fig.~\ref{Fig:
dis fact}.
\end{theorem}
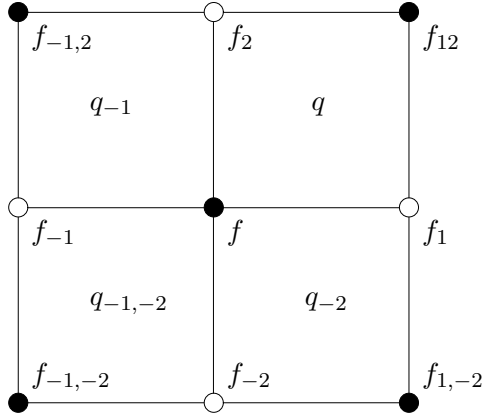
\begin{figure}[htbp]
\begin{center}
\setlength{\unitlength}{0.045em}
\begin{picture}(300,320)(0,0)
 \put(0,0){\circle*{15}} \put(150,0){\circle{15}}
 \put(300,0){\circle*{15}}
 \put(0,150){\circle{15}} \put(150,150){\circle*{15}}
 \put(300,150){\circle{15}}
 \put(0,300){\circle*{15}} \put(150,300){\circle{15}}
 \put(300,300){\circle*{15}}
 \path(0,0)(142.5,0) \path(157.5,0)(300,0)
 \path(0,0)(0,142.5) \path(0,157.5)(0,300)
 \path(150,7.5)(150,292.5)
 \path(7.5,150)(292.5,150)
 \path(0,300)(142.5,300)    \path(300,0)(300,142.5)
 \path(157.5,300)(300,300)
 \path(300,157.5)(300,300)
 \put(160,125){$f$}
 \put(310,125){$f_1$}
 \put(10,125){$f_{-1}$}
 \put(160,15){$f_{-2}$}
 \put(310,15){$f_{1,-2}$}
 \put(10,15){$f_{-1,-2}$}
 \put(160,275){$f_2$}
 \put(310,275){$f_{12}$}
 \put(10,275){$f_{-1,2}$}
 \put(55,75){$q_{-1,-2}$}
 \put(55,225){$q_{-1}$}
 \put(220,75){$q_{-2}$}
 \put(225,225){$q$}
\end{picture}
\caption{Four adjacent quadrilaterals of a discrete isothermic
net: the cross-ratios satisfy $q\cdot q_{-1,-2}=q_{-1}\cdot
q_{-2}$, the five points $f$ and $f_{\pm 1,\pm 2}$ lie on a
sphere}\label{Fig: dis fact}
\end{center}
\end{figure}
{\bf Proof.} Perform a M\"obius transformation sending $f$ to
$\infty$. Under such a transformation, the four adjacent circles
through $f$ turn into four straight lines $f_{\pm 1}f_{\pm 2}$,
containing the corresponding points $f_{\pm 1,\pm 2}$. Formula
(\ref{eq:i/j prod}) turns into the following relation for the
quotients of (directed) lengths:
\begin{equation}\label{eq: 8-ratio}
\frac{|f_2f_{12}|}{|f_{12}f_1|}\cdot
\frac{|f_1f_{1,-2}|}{|f_{1,-2}f_{-2}|}\cdot
\frac{|f_{-2}f_{-1,-2}|}{|f_{-1,-2}f_{-1}|}\cdot
\frac{|f_{-1}f_{-1,2}|}{|f_{-1,2}f_2|}=1.
\end{equation}
If the affine space through the points $f_{\pm 1}$, $f_{\pm 2}$ is
three-dimesnional, then eq. (\ref{eq: 8-ratio}) is equivalent to
the fact that the four points $f_{\pm 1,\pm 2}$ lie in a plane,
which is a sphere through $f=\infty$. This is a three-dimensional
generalization of the Menelaus theorem; since this generalization
is not very well known, we give it, with a proof, in appendix to
this section. If, on the contrary, the four points $f_{\pm 1}$,
$f_{\pm 2}$ are co-planar, then we are in the situation of
Fig.~\ref{Fig: Desargues}, described by the Desargues theorem.
Here, we apply the Menelaus theorem twice, to the triangle
$f_{-1}f_2f_1$ intersected by the line $f_{-1,2}f_{12}$, and to
the triangle $f_{-1}f_{-2}f_1$ intersected by the line
$f_{-1,-2}f_{1,-2}$: both lines meet the line $f_{-1}f_1$ at the
same point $\ell^{(1)}$, if and only if
\[
\frac{|f_2f_{12}|}{|f_{12}f_1|}\cdot
\frac{|f_{-1}f_{-1,2}|}{|f_{-1,2}f_2|}
=-\frac{|f_{-1}\ell^{(1)}|}{|\ell^{(1)}f_1|}=
\frac{|f_{-2}f_{1,-2}|}{|f_{1,-2}f_1|}\cdot
\frac{|f_{-1}f_{-1,-2}|}{|f_{-1,-2}f_{-2}|}.
\]
This yields (\ref{eq: 8-ratio}). $\Box$
\smallskip

The claim of Theorem \ref{thm: isothermic cross-ratios} can be
re-formulated as follows.

\begin{corollary}\label{cor: isothermic cross-ratios}
{\bf (Factorized cross-ratios for a discrete isothermic net)} A
circular net $f:\bbZ^2\to\bbR^N$ is a discrete isothermic net, if
and only if the cross-ratios $q=q(f,f_1,f_{12},f_2)$ of its
elementary quadrilaterals satisfy the following condition:
\begin{equation}\label{eq:1/2}
q(f,f_1,f_{12},f_2)=\frac{\alpha_1}{\alpha_2}\,,
\end{equation}
with some edge functions $\alpha_i$ satisfying the labelling
property (\ref{eq: labelling}).
\end{corollary}
Clearly, functions $\alpha_i$ are defined up to a common constant
factor. Actually, it was this characterization of discrete
isothermic nets that was used as a definition in the pioneering
paper \cite{BP1}.

Actually, edge functions $\alpha_i$ in Corollary \ref{cor:
isothermic cross-ratios} admit a nice geometric expression.
According to Theorem \ref{Thm: charact 2d T-net}, a discrete
isothermic net $f:\bbZ^2\to\bbR^N$ can be characterized by the
existence of representatives $\hat{y}:\bbZ^2\to\bbL^{N+1,1}$ in
the light cone satisfying the discrete Moutard equation
(\ref{eq:dMou 2d minus}). We will use for these T-net
representatives the notation
\begin{equation}\label{eq: T-net in L repr}
\hat{y}=s^{-1}\hat{f}=s^{-1}(f+\ee_0+|f|^2\ee_\infty).
\end{equation}
Thus, $s^{-1}$ denotes the $\ee_0$-component of the T-net
representative $\hat{y}$ of an isothermic net $f$. The function
$s:\bbZ^2\to\bbR$ plays the role of a discrete metric of the
isothermic net $f$. Now, eq. (\ref{eq: alpha dMou quadr}) defines
functions on edges,
\begin{equation}\label{eq: iso Mou labelling}
\alpha_i=-2\langle
\hat{y},\hat{y}_i\rangle=\frac{|f_i-f|^2}{ss_i},
\end{equation}
possessing property (\ref{eq: labelling}): any two opposite edges
of any elementary square carry the same value of the corresponding
$\alpha_i$.

\begin{theorem}
{\bf (Cross-ratios through discrete metric)} Edge functions
$\alpha_i$ participating in the factorization (\ref{eq:1/2}) of
the cross-ratios of elementary quadrilaterals of a discrete
isothermic net can be defined by eq. (\ref{eq: iso Mou
labelling}).
\end{theorem}
{\bf Proof.} Comparing the $\ee_0$-components in the Moutard
equation $\hat{y}_{12}-\hat{y}=a_{12}(\hat{y}_2-\hat{y}_1)$, we
find: $a_{12}=(s_{12}^{-1}-s^{-1})/(s_2^{-1}-s_1^{-1})$.
Therefore, we can re-write the Moutard equation as
\[
\Big(\frac{1}{s_2}-\frac{1}{s_1}\Big)\Big(\frac{\hat{f}_{12}}{s_{12}}
-\frac{\hat{f}}{s}\Big)=
\Big(\frac{1}{s_{12}}-\frac{1}{s}\Big)\Big(\frac{\hat{f}_2}{s_2}
-\frac{\hat{f}_1}{s_1}\Big),
\]
which is equivalent to
\begin{equation}\label{eq: Mou in cone}
\frac{\hat{f}_1-\hat{f}}{ss_1}+
\frac{\hat{f}_{12}-\hat{f}_1}{s_1s_{12}}=\frac{\hat{f}_2-\hat{f}}{ss_2}+
\frac{\hat{f}_{12}-\hat{f}_2}{s_2s_{12}}\,.
\end{equation}
The $\bbR^N$-part  of the latter equation, i.e.,
\begin{equation}\label{eq: Mou in cone 1}
\frac{f_1-f}{ss_1}+\frac{f_{12}-f_1}{s_1s_{12}}=
\frac{f_2-f}{ss_2}+ \frac{f_{12}-f_2}{s_2s_{12}}\,,
\end{equation}
can be rewritten with the help of eq. (\ref{eq: iso Mou
labelling}) as
\begin{equation}\label{eq: Mou in cone 2}
\alpha_1\frac{f_1-f}{|f_1-f|^2}+\alpha_2\frac{f_{12}-f_1}{|f_{12}-f_1|^2}=
\alpha_2\frac{f_2-f}{|f_2-f|^2}+\alpha_1\frac{f_{12}-f_2}{|f_{12}-f_2|^2}.
\end{equation}
In terms of the inversion in the Clifford algebra
$\cC\ell(\bbR^N)$, this can be presented as
\begin{equation}\label{eq: Mou in cone 3}
\alpha_1(f_1-f)^{-1}+\alpha_2(f_{12}-f_1)^{-1}=
\alpha_2(f_2-f)^{-1}+\alpha_1(f_{12}-f_2)^{-1}.
\end{equation}
This latter equation is, in the generic case $f_{12}+f\neq
f_1+f_2$, equivalent to eq. (\ref{eq:1/2}). It is not quite
straightforward to show this equivalence in case of
non-commutative variables $f\in\cC\ell(\bbR^N)$. But one can
identify the plane of the quadrilateral $(f,f_1,f_{12},f_2)$ with
$\bbC$, and then eq. (\ref{eq: Mou in cone 2}) is the (complex
conjugate of) eq. (\ref{eq: Mou in cone 3}), where now all
variables are commutative (complex numbers), and in this case the
equivalence to eq. (\ref{eq:1/2}) is immediate. $\Box$
\smallskip

T-nets in the light cone $\bbL^{N+1,1}$ are 3D-consistent. This
yields also the 3D-consistency of the cross-ratio equation
(\ref{eq:1/2}) with prescribed labelling $\alpha_i$ of the edges,
i.e., of the 2D equation
\begin{equation}\label{eq:i/j}
q(f,f_i,f_{ij},f_j)=\frac{\alpha_i}{\alpha_j}.
\end{equation}
Both constructions provide us with a well-defined notion of
multidimensional discrete isothermic nets, and therefore with
Darboux transformations of discrete isothermic nets with the usual
permutability properties.
\smallskip

We finish this section with the notion of duality for discrete
isothermic nets.
\begin{theorem}\label{Th: dis dual}
{\bf (Dual discrete isothermic net)} Let $f:\bbZ^m\to\bbR^N$ be a
discrete isothermic net, with the T-net representatives in the
light cone
\[
\hat{y}=s^{-1}\hat{f}=s^{-1}(f+\ee_0+|f|^2\ee_\infty):\,\bbZ^m\to\bbL^{N+1,1}.
\]
Then the $\bbR^N$-valued discrete one-form $\delta f^*$ defined by
\begin{equation}\label{eq: dis dual ij}
\delta_i f^*=\frac{\delta_i f}{ss_i}=\alpha_i\frac{\delta_i
f}{|\delta_i f|^2}\,, \qquad i=1,\ldots,m,
\end{equation}
is closed. Its integration defines (up to translation) a net
$f^*:\,\bbZ^m\to\bbR^N$, called dual to the net $f$. The net $f^*$
is a discrete isothermic net, with
\begin{equation}\label{eq:dis prop dual}
 q(f^*,f_i^*,f_{ij}^*,f_j^*)=\frac{\alpha_i}{\alpha_j}\,.
\end{equation}
Define also the function $s^*:\,\bbZ^m\to\bbR$ as $s^*=s^{-1}$.
Then the net
\[
\hat{y}^*=(s^*)^{-1}\hat{f^*}=
(s^*)^{-1}(f^*+\ee_0+|f^*|^2\ee_\infty): \,\bbZ^m\to\bbL^{N+1,1}
\]
is a T-net in the light cone.
\end{theorem}
{\bf Proof.} Clearly, for any pair of indices $i,j$ the function
$\hat{f}$ satisfies an equation analogous to eq. (\ref{eq: Mou in
cone}), which expresses the closeness of the $\bbR^{N+1,1}$-valued
one-form defined by $\delta_i\hat{f}/(ss_i)$. Unfortunately, the
net obtained by integration of this one-form does not lie, in
general, in the light cone $\bbL^{N+1,1}$ and cannot be taken as
the dual net $\hat{f}^*$. We use the following trick for the
construction of the dual net $\hat{f}^*$ in the light cone. The
$\bbR^N$-part of eq. (\ref{eq: Mou in cone}), i.e., eq. (\ref{eq:
Mou in cone 1}), expresses the closeness of the $\bbR^{N}$-valued
one-form $\delta_if^*=\delta_if/(ss_i)$, whose integration gives a
dual net $f^*$ in $\bbR^N$. Eq. (\ref{eq:dis prop dual}) follows
immediately from eq. (\ref{eq: dis dual ij}) and implies that
$f^*$ is a discrete isothermic net. In particular, it is a
circular net, so that $\hat{f}^*=f^*+\ee_0+|f^*|^2\ee_\infty$ is a
conjugate net in the light cone. It remains to show that the so
defined $\hat{f}^*$ is a Moutard net, with the T-net
representatives $\hat{y}^*=(s^*)^{-1}\hat{f}^*$. This claim is
equivalent to closeness of the discrete $\bbR^{N+1,1}$-valued
one-form $\delta_i\hat{f}^*/(s^*s_i^*)$. Since $\hat{f}^*$ is a
conjugate net in the light cone, it is enough to prove the
closeness of the $\bbR^N$-valued one-form
$\delta_if^*/(s^*s_i^*)$. But for $s^*=s^{-1}$ we have:
\[
\frac{\delta_if^*}{s^*s_i^*}=(ss_i)\delta_if^*=
(ss_i)\frac{\delta_if}{ss_i}=\delta_if,
\]
which is automatically closed. $\Box$

\subsection{Appendix: generalized Menelaus theorem}

\begin{theorem}\label{prop: gen Menelaus}
Let $P_0$, ..., $P_n$ be $n+1$ points in general position in
$\bbR^{n}$, so that the affine space through the points $P_i$ is
$n$-dimensional. Let $P_{i,i+1}$ be some points on the lines
$P_iP_{i+1}$ (indices are read modulo $n+1$). The $n+1$ points
$P_{i,i+1}$ lie in an $(n-1)$-dimensional affine subspace, if and
only if the following relation for the quotients of the directed
lengths holds:
\[
\prod_{i=0}^n
\frac{|P_iP_{i,i+1}|}{|P_{i,i+1}P_{i+1}|}=(-1)^{n+1}.
\]
\end{theorem}
{\bf Proof.} The points $P_{i,i+1}$ lie in an $(n-1)$-dimensional
affine subspace, if there is a non-trivial linear dependence
\[
\sum_{i=0}^n\alpha_iP_{i,i+1}=0 \quad {\rm with} \quad
\sum_{i=0}^n\alpha_i=0.
\]
Substituting $P_{i,i+1}=(1-\xi_i)P_i+\xi_iP_{i+1}$, and taking
into account the general position condition, which can be read as
linear independence of the vectors $P_i-P_0$, we come to a
homogeneous system of $n+1$ linear equations for $n+1$
coefficients $\alpha_i$:
\[
\xi_i\alpha_i+(1-\xi_{i+1})\alpha_{i+1}=0,\quad i=0,\ldots,n
\]
(where indices are understood modulo $n+1$). Clearly it admits a
non-trivial solution if and only if
\[
\prod_{i=0}^n
\frac{\xi_i}{1-\xi_i}=\prod_{i=0}^n\frac{|P_iP_{i,i+1}|}{|P_{i,i+1}P_{i+1}|}
=(-1)^{n+1}.\qquad \Box
\]

\section{Isothermic surfaces in Laguerre geometry}
\label{subsect: Lag isothermic}

Definitions and constructions of this section can be generalized
for higher-dimensional nets
\[
P:\bbZ^m\to\{{\rm hyperplanes\ in\ }\bbR^N\},
\]
with the most natural geometric case $m=N-1$ (discrete Laguerre
geometry of hypersurfaces in $\bbR^N$); we restrict ourselves to
the case of surfaces in $\bbR^3$, i.e., $N=3$, $m=2$.

\begin{definition}\label{def: Lag isothermic}
{\bf (Discrete L-isothermic surface)} A two-dimensional conical
net $P:\bbZ^2\to\{\rm planes\ in\ \bbR^3\}$ is called a discrete
L-isothermic surface, if the corresponding
$\hat{p}:\bbZ^2\to\bbL^{4,2}$ is a discrete Moutard net.
\end{definition}
Recall that for an (oriented) plane $P=\{x\in\bbR^3: \langle
v,x\rangle=d\}$ with the unit normal vector $v\in\bbS^2$ and
$d\in\bbR$ its representative $\hat{p}\,$ in the Lie quadric
$\bbL^{4,2}$ is given by
\[
\hat{p}=v+0\cdot\ee_0+2d\ee_\infty+1\cdot\ee_6.
\]
Recall also that the vectors $v:\bbZ^2\to\bbS^2$ comprise the {\em
Gauss map} for a given conical net $P:\bbZ^2\to\{\rm planes\ in\
\bbR^3\}$.

From Definitions \ref{def: dmn}, \ref{def: dmn 3d} there follows a
geometric characterization of discrete L-isothermic nets:
\begin{theorem}\label{thm: L-isothermic 5-plane sphere}
{\bf (Central spheres for discrete L-isothermic nets)}

(a) A conical net $P:\bbZ^2\to\{\rm planes\ in\ \bbR^3\}$ not
tangent to a two-dimensional sphere is a discrete L-isothermic
net, if and only if for any $u\in\bbZ^2$ the five planes $P$ and
$P_{\pm 1,\pm 2}$ are tangent to a two-dimensional sphere not
touching some (and then any) of the four planes $P_{\pm 1}$,
$P_{\pm 2}$.

(b) A conical net $P:\bbZ^2\to\{\rm tangent\ planes\ of\
\bbS^2\subset\bbR^3\}$ is a discrete L-isothermic net, if and only
if for any $u\in\bbZ^2$ the three cones through $P$,
\begin{eqnarray*}
& C^{\rm(up)}={\rm cone}(P,P_{12},P_{-1,2}),\quad
C^{\rm(down)}={\rm cone}(P,P_{1,-2},P_{-1,-2}), &
\\
& C^{\rm(1)}={\rm cone}(P,P_{1},P_{-1}), &
\end{eqnarray*}
have one additional plane in common, which is also equivalent for
the three cones through $P$,
\begin{eqnarray*}
& C^{\rm(left)}={\rm cone}(P,P_{-1,2},P_{-1,-2}),\quad
C^{\rm(right)}={\rm cone}(P,P_{1,2},P_{1,-2}), &
\\
& C^{\rm(2)}={\rm cone}(P,P_{2},P_{-2}), &
\end{eqnarray*}
to have one additional plane in common.
\end{theorem}

\begin{figure}[htbp]
\begin{center}
{\includegraphics[width=0.7\textwidth]{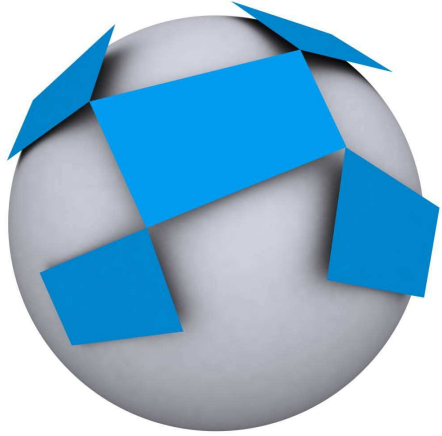}}
\end{center}
\caption{Five diagonally neighboring planes of a generic discrete
L-isothermic surface, with a central sphere.}
\label{fig:L-isothermic}
\end{figure}

The (generic) case a) of Theorem \ref{thm: L-isothermic 5-plane
sphere} is illustrated on Fig. \ref{fig:L-isothermic}.

\begin{theorem}
{\bf (Gauss map of an L-isothermic net is an isothermic net in the
sphere)} Gauss map of an L-isothermic net is a discrete isothermic
net in $\bbS^2$. Conversely, if for any $u\in\bbZ^2$ the four
planes $P,P_1,P_2,P_{12}$ of a net $P:\bbZ^2\to\{\rm planes\ in\
\bbR^3\}$ meet at a point, and the Gauss map of the net $P$ is
isothermic, then $P$ is an L-isothermic conical net.
\end{theorem}
{\bf Proof.} First, let $P$ be an L-isothermic net. Then for some
$c:\bbZ^2\to\bbR$ the net $c\hat{p}$ is a T-net in the Lie
quadric. As a consequence, $c(v+1\cdot\ee_0+1\cdot\ee_\infty)$ is
a T-net in the light cone $\bbL^{4,1}$ of the Minkowski space
$\bbR^{4,1}$ of the M\"obius geometry for $N=3$. Therefore, the
net $v:\bbZ^2\to\bbS^2\subset\bbR^3$ is isothermic.

Conversely, let the net $v:\bbZ^2\to\bbS^2\subset\bbR^3$ be
isothermic. This is equivalent to the existence of the function
$c:\bbZ^2\to\bbR$ such that $c(v,1)$ is a T-net. If now $\langle
v,x\rangle=d$ is the equation of the plane $P$, then the existence
of the common intersection point of the planes $P,P_1,P_2,P_{12}$
yields that the function $cd$ satisfies the same Moutard equation
as the function $cv$. Therefore, $c(v,d,1)$ is a T-net, so that
$\hat{p}\,$ is a discrete Moutard net. $\Box$

\section{Lie geometry: S-isothermic nets}

Two-dimensional nets in the Lie quadric $\bbL^{4,2}$ are discrete
congruences of spheres. An interesting class of such congruences
is constituted by discrete Moutard nets in $\bbP(\bbL^{4,2})$. We
leave a general study of this class for a future research, and
describe here as an example a particularly interesting subclass,
for which the T-net representatives in $\bbL^{4,2}$ have a fixed
$\ee_6$-component:
\[
\hat{y}=\frac{\kappa}{r}\Big(c+\ee_0+(|c|^2-r^2)\ee_\infty+r\ee_6\Big).
\]
Omitting the constant and therefore non-interesting
$\ee_6$-component, we come to a T-net in a hyperboloid of the
Lorentz space of the M\"obius geometry,
\[
\bbL^{4,1}_\kappa=\big\{\xi\in\bbR^{4,1}:\langle
\xi,\xi\rangle=\kappa^2\big\}.
\]

\begin{definition}\label{dfn:S-iso}
{\bf (S-isothermic net)} A map
\[
S:\bbZ^2\to\{{\mbox {\rm oriented spheres in\;\,}} \bbR^3\}
\]
is called an S-isothermic net, if the corresponding map
\begin{equation}\label{eq: s in kappa}
\hat{s}:\bbZ^2\to\bbL^{4,1}_\kappa,\qquad
\hat{s}=\frac{\kappa}{r}\Big(c+\ee_0+(|c|^2-r^2)\ee_\infty\Big),
\end{equation}
is a T-net.
\end{definition}

Thus, S-isothermic nets are governed by equation
\begin{equation}\label{eq: S-iso Mou a}
\hat{s}_{12}-\hat{s}=a_{12}(\hat{s}_2-\hat{s}_1),\qquad
a_{12}=\frac{\langle
\hat{s},\hat{s}_1-\hat{s}_2\rangle}
{\kappa^2-\langle\hat{s}_1,\hat{s}_2\rangle}=
\frac{\alpha_1-\alpha_2}{\kappa^2-\langle
\hat{s}_1,\hat{s}_2\rangle}\,,
\end{equation}
with the quantities $\alpha_i=\langle \hat{s},\hat{s}_i\rangle$
depending on $u_i$ only. If (oriented) radii of all hyperspheres
become uniformly small, $r(u)\sim\kappa s(u)$, $\kappa\to 0$, then
in the limit we recover discrete isothermic nets.

Consistency of T-nets in $\bbL^{4,1}_{\,\kappa}$ (which is a
particular case of Theorem \ref{Th: M-lat in quadric}) yields, in
particular, Darboux transformations for S-isothermic nets. A
Darboux transform $\hat{s}^+:\bbZ^m\to\bbL^{4,1}_{\,\kappa}$ of a
given S-isothermic net $\hat{s}:\bbZ^m\to\bbL^{4,1}_{\,\kappa}$ is
uniquely specified by a choice of one of its spheres
$\hat{s}^+(0)$.

We turn now to geometric properties of S-isothermic nets. First of
all, S-isothermic nets form a subclass of discrete R-congruences
of spheres  (see \cite{BS2} for a geometric characterization of
discrete R-congruences). Further, consider the quantities
$\langle\hat{s},\hat{s}_i\rangle$ which have the meaning of
cosines of the intersection angles of the neighboring spheres
(resp., of their so called inversive distances if they do not
intersect). Then these quantities
$\langle\hat{s},\hat{s}_i\rangle$ have the labelling propery,
i.e., depend on $u_i$ only.

There holds the following generalization of Theorem \ref{Th: dis
dual}.
\begin{theorem}\label{Th: S-iso dual}
{\bf (Dual S-isothermic net)} Let
\[
S:\bbZ^m\to\{{\mbox {\rm oriented spheres in\ }}\bbR^3\}
\]
be an S-isothermic net. Denote the Euclidean centers and
(oriented) radii of $S$ by $c:\bbZ^m\to\bbR^3$ and
$r:\bbZ^m\to\bbR$, respectively. Then the $\bbR^3$-valued discrete
one-form $\delta c^*$ defined by
\begin{equation}\label{eq: S-iso dual}
\delta_i c^*=\frac{\delta_i c}{rr_i}\,,\qquad 1\le i\le m,
\end{equation}
is closed, so that its integration defines (up to a translation) a
function $c^*:\bbZ^m\to\bbR^3$. Define also $r^*:\bbZ^m\to\bbR$ by
$r^*=r^{-1}$. Then the spheres $S^*$ with the centers $c^*$ and
radii $r^*$ form an S-isothermic net, called {\em dual} to $S$.
\end{theorem}
{\bf Proof.} Consider equation
\begin{equation}\label{eq:minusMou}
\hat{s}_{ij}-\hat{s}=a_{ij}(\hat{s}_j-\hat{s}_i),
\end{equation}
in terms of $\hat{s}$ from (\ref{eq: s in kappa}). Its
$\ee_0$-part yields:
$a_{ij}=(r_{ij}^{-1}-r^{-1})/(r_j^{-1}-r_i^{-1})$. This allows us
to rewrite eq. (\ref{eq:minusMou}) as
\begin{equation}\label{eq: S-iso dual eq}
(r_j^{-1}-r_i^{-1})(\hat{s}_{ij}-\hat{s})=
(r_{ij}^{-1}-r^{-1})(\hat{s}_j-\hat{s}_i).
\end{equation}
A direct computation shows that the $\bbR^3$-part of this equation
can be rewritten as
\begin{equation}
\frac{c_i-c}{rr_i}+\frac{c_{ij}-c_i}{r_ir_{ij}}=
\frac{c_j-c}{rr_j}+\frac{c_{ij}-c_j}{r_jr_{ij}}\,,
\end{equation}
which is equivalent to closeness of the form $\delta c^*$ defined
by (\ref{eq: S-iso dual}). In the same way, the $\ee_\infty$-part
of eq. (\ref{eq: S-iso dual eq}) is equivalent to closeness of the
discrete form $\delta w$ defined by
\begin{equation*}\label{eq: S-iso dual w}
\delta_i w=\frac{\delta_i(|c|^2-r^2)}{rr_i}\,,\qquad 1\le i\le m.
\end{equation*}
For similar reasons, the second claim of the theorem is equivalent
to closeness of the form
\begin{equation*}\label{eq: S-iso dual z}
\delta_i w^*=\frac{\delta_i\big(|c^*|^2-(r^*)^2\big)}{r^*r_i^*}\,,
\qquad 1\le i\le m,
\end{equation*}
where, recall, $r^*=1/r$. With the help of
$c_i^*-c^*=(c_i-c)/rr_i$, one easily checks that the forms $\delta
w$ and $\delta w^*$ can be written as
\begin{eqnarray*}
\delta_i w & = & \langle
c_i^*-c^*,c_i+c\rangle-\frac{r_i}{r}+\frac{r}{r_i}\,,\\
\delta_i w^* & = & \langle
c_i-c,c_i^*+c^*\rangle-\frac{r}{r_i}+\frac{r_i}{r}\,.
\end{eqnarray*}
The sum of these one-forms is closed:
\[
\delta_i(w+w^*)=2\langle c_i^*,c_i\rangle-2\langle c^*,c\rangle,
\]
therefore they are closed simultaneously. $\Box$
\medskip

An interesting particular case of S-isothermic surfaces is
characterized by touching of any pair of neighboring spheres. In
this case the limit of small spheres is not relevant, therefore it
is convenient to restrict the considerations to a fixed value of
$\kappa=1$. Clearly, in this case both $\alpha_i=\langle
\hat{s},\tau_i\hat{s}\rangle$, $i=1,2$, can, in principle, take
values $\pm 1$. However, it is easily seen from (\ref{eq: S-iso
Mou a}) that in case $\alpha_1=\alpha_2$ one gets only trivial
nets. Thus, we assume that
\begin{equation}\label{eq: S-iso touch}
\langle\hat{s},\hat{s}_1\rangle=\langle\hat{s}_2,\hat{s}_{12}\rangle
=-1,\quad \langle\hat{s},\hat{s}_2\rangle
=\langle\hat{s}_1,\hat{s}_{12}\rangle=1.
\end{equation}
Interestingly, these touching conditions already yield that the
linear dependence of the spheres has the Moutard shape.
\begin{theorem}
S-isothermic surfaces with touching spheres can be characterized
by any of the following equivalent descriptions:
\begin{itemize}
    \item Q-congruence of spheres (Q-net in the Lorentz space of M\"obius geometry)
    with touching spheres
    \item R-congruence (Q-net in Lie quadric) with touching
    spheres
    \item T-net in the Lie quadric with touching spheres
\end{itemize}
\end{theorem}
{\bf Proof.} Let $\hat{s}$, $\hat{s}_1$, $\hat{s}_2$,
$\hat{s}_{12}$ be four oriented spheres in $\bbL^{4,1}_1$,
pairwise touching so that eq. (\ref{eq: S-iso touch}) is
fulfilled, and linearly dependent. (We remark that in the present
situation the geometric meaning of linear dependence is the
existence of a common orthogonal circle through the touching
points.) We show that the linear dependence has to be of the
Moutard form:
\begin{equation}\label{eq: S-iso touch eq}
\hat{s}_{12}-\hat{s}=a_{12}(\hat{s}_2-\hat{s}_1),\qquad
a_{12}=-\frac{2}{1-\langle\hat{s}_1,\hat{s}_2\rangle}.
\end{equation}
In this proof, we make a general position assumption that the
spheres $\hat{s}$ and $\hat{s}_{12}$ do not touch and the spheres
$\hat{s}_1$ and $\hat{s}_2$ do not touch. Write the linear
dependence condition as
\begin{equation}\label{eq: S-iso aux}
\hat{s}_{12}=\lambda\hat{s}+\mu \hat{s}_1+\nu\hat{s}_2.
\end{equation}
Scalar product of this with $\hat{s}_1$, $\hat{s}_2$ leads to:
\[
1+\lambda=\mu+\nu\langle\hat{s}_1,\hat{s}_2\rangle
=-\mu\langle\hat{s}_1,\hat{s}_2\rangle-\nu
 \quad\Rightarrow\quad \mu=-\nu=\frac{\lambda+1}
{1-\langle\hat{s}_1,\hat{s}_2\rangle}
\]
Similarly, a scalar product of eq. (\ref{eq: S-iso aux}) with
$\hat{s}$, $\hat{s}_{12}$ leads to:
\[
\mu-\nu=\lambda-\langle\hat{s},\hat{s}_{12}\rangle=
1-\lambda\langle\hat{s},\hat{s}_{12}\rangle \quad\Rightarrow\quad
\lambda=1.
\]
This yields eq. (\ref{eq: S-iso touch eq}). $\Box$

\section{Remarks on the continuous limit}

To perform a continuous limit in constructions of Discrete
Differential Geometry, one should think of the underlying lattice
$\bbZ^2$ of discrete surfaces as of $(\epsilon\bbZ)^2$ and then
send $\epsilon\to 0$. In such a limit it is common to assume that
the discrete functions $y(n_1,n_2)$ on $\bbZ^2$ approximate
sufficiently smooth functions $y(u_1,u_2)$ on $\bbR^2$, if
$n_i\epsilon=u_i$.

{\bf Moutard nets.} It turns out, however, that the discrete
Moutard equation in the form (\ref{eq:dMou 2d minus}) with minus
signs does not admit a nice continuous limit, while its close
relative -- discrete Moutard equation with plus signs -- does:
\begin{equation}\label{eq:dMou 2d plus}
y_{12}+y=a_{12} (y_1+y_2).
\end{equation}
Indeed, if $a_{12}=1+\frac{\epsilon^2}{4}q_{12}$, then the above
equation is re-written as
\[
\frac{1}{\epsilon^2}(y_{12}-y_1-y_2+y)=\frac{1}{4}q_{12}(y_{12}+y_1+y_2+y),
\]
which in the limit $\epsilon\to 0$ tends to the {\em Moutard
equation} \cite{M}:
\begin{equation}\label{eq: Mou}
\partial_1\partial_2 y=qy.
\end{equation}
Equations (\ref{eq:dMou 2d minus}) and (\ref{eq:dMou 2d plus}) are
actually related by a very simple transformation of dependent
variables: $y(n_1,n_2)\mapsto(-1)^{n_2}y(n_1,n_2)$. As long as $y$
play the role of homogeneous coordinates, such a transformation
does not change the corresponding points in the projective space,
therefore in the formulation of Theorem \ref{Thm: charact 2d
T-net} we could use eq. (\ref{eq:dMou 2d plus}) in place of eq.
(\ref{eq:dMou 2d minus}), as well. However, eq. (\ref{eq:dMou 2d
plus}) cannot be posed on all faces of a three-dimensional
lattice, and it is not multidimensionally consistent. As a
consequence, the right way to construct discrete Moutard
transformations for nets governed by eq. (\ref{eq:dMou 2d plus})
is to consider multidimensional T-nets and then perform a change
of variables $y(n)\mapsto(-1)^{n_2}y(n)$ which {\em breaks the
symmetry} among lattice directions. This leads to changing the
signs in some of the formulas (\ref{eq:dmp T}) for discrete
Moutard transformations, which turn into
\begin{equation}\label{eq:dmp}
y^+_1-y=b_1(y^+-y_1),\qquad  y^+_2+y=b_2(y^++y_2).
\end{equation}
An advantage of this break of symmetry is again the possibility to
perform a continuous limit in eqs. (\ref{eq:dmp}), recovering the
formulas of the {\em Moutard transformation} for eq. (\ref{eq:
Mou}), see \cite{M}:
\begin{equation}\label{eq:mp}
\partial_1 y^++\partial_1 y=p_1(y^+-y),\qquad \partial_2 y^+-\partial_2 y
=p_2(y^++y).
\end{equation}

{\bf Isothermic nets.} A similar situation takes place for
discrete isothermic nets. In order to enable the continuous limit
to smooth isothermic surfaces, one should start with a discrete
two-dimensional isothermic net with embedded elementary
quadrilaterals. It is convenient to represent their negative
cross-ratios as
\begin{equation}\label{eq:dis property}
  q(f,f_1,f_{12},f_2)=-\frac{\alpha_1}{\alpha_2}\,,
\end{equation}
with positive labels $\alpha_1$ and $\alpha_2$. Formally, this
means nothing more than changing the notation $\alpha_2\mapsto
-\alpha_2$. But this operation puts some formulas into the form
which allows for a continuous limit. If one keeps the both
formulas
\[
|f_1-f|^2=\alpha_1ss_1, \qquad |f_2-f|^2=\alpha_2ss_2,
\]
then this implies a slight modification in the definition of the
function $s$, namely $s(u)\mapsto (-1)^{u_2}s(u)$, which assures
the positivity of $s$ in case of positive $\alpha_1$, $\alpha_2$.
Only upon this modification does the positive function $s$ have a
continuous limit. In this limit we recover a characterization of
isothermic surfaces as those curvature line parametrized surfaces
for which the first fundamental form is conformal, possibly upon a
re-parametrization of curvature lines:
\begin{equation}\label{eq:is prop}
\langle\partial_1 f,\partial_2
f\rangle=0,\quad |\partial_1 f|^2=\alpha_1 s^2,\quad |\partial_2
f|^2=\alpha_2 s^2,
\end{equation}
with $\alpha_i=\alpha_i(u_i)$. Eq. (\ref{eq: dis dual ij}) for the
dual surface turns into
\begin{equation}\label{eq: dis dual}
\delta_1 f^*=\alpha_1\frac{\delta_1 f}{|\delta_1 f|^2}=
\frac{\delta_1 f}{ss_1}\,,  \qquad \delta_2
f^*=-\alpha_2\frac{\delta_2 f}{|\delta_2 f|^2}= -\frac{\delta_2
f}{ss_2}\,,
\end{equation}
which is a direct discrete analogue of equations
\begin{equation}\label{eq: is dual}
\partial_1 f^*=\alpha_1\frac{\partial_1 f}{|\partial_1 f|^2}=
\frac{\partial_1 f}{s^2},\qquad
\partial_2 f^*=-\alpha_2\frac{\partial_2 f}{|\partial_2 f|^2}=
-\frac{\partial_2 f}{s^2},
\end{equation}
defining the dual isothermic surfaces in the smooth case.

{\bf Remark.} In the smooth case the functions $\alpha_1$,
$\alpha_2$ in eq. (\ref{eq:is prop}) can be absorbed into a
re-parametrization of the independent variables
$u_i\mapsto\varphi_i(u_i)$ $\,(i=1,2)$, by which one can always
achieve that $\alpha_1=\alpha_2=1$, so that the first fundamental
form of the surface $f$ is conformal. Of course, such a
re-parametrization is not possible in the discrete context.
Nevertheless, one can consider a narrower class of discrete
isothermic surfaces, characterized by eq. (\ref{eq:dis property})
with $\alpha_1=\alpha_2=1$:
\begin{equation}\label{eq: cr=-1}
 q(f,f_1,f_{12},f_2)=-1.
\end{equation}
This condition (all elementary quadrilaterals of $f$ are conformal
squares) may be regarded as a discretization of the conformality
of the first fundamental form. Eq. (\ref{eq: cr=-1}), being a
particular case of eq. (\ref{eq:i/j}) with a special labelling,
enjoys all the properties of the general case. However, it is
important to observe that it is not 3D consistent {\em with
itself}, i.e., it cannot be imposed on all faces of a 3D cube.
Indeed, if $\alpha_1/\alpha_2=-1$, then it is impossible to have
additionally $\alpha_2/\alpha_3=-1$ and $\alpha_1/\alpha_3=-1$.
\smallskip

Note that the above change of signs yields also the similar
modification in the lift $\hat{s}$, which therefore satisfies the
discrete Moutard equation with plus signs (\ref{eq:dMou 2d plus}),
admitting the continuous limit. This characterization of
isothermic surfaces is due to Darboux \cite{Da}:
\begin{theorem}\label{thm: Darboux}
{\bf (Isothermic surfaces as Moutard nets in the light cone)} A
surface $f:\bbR^2\to\bbR^N$ is isothermic, if and only if there
exists a function $s:\bbR^2\to\bbR$ such that the lift of $f$ into
the light cone
\[
\hat{s}=s^{-1}(f+\ee_0+|f|^2\ee_\infty):\bbR^2\to\bbL^{N+1,1}
\]
satisfies a Moutard equation (\ref{eq: Mou}). Darboux
transformations of an isothermic surface $f$ are obtained by
projection to $\bbR^N$ from Moutard transformations of $\hat{s}$
within the light cone $\bbL^{N+1,1}$.
\end{theorem}

{\bf L-isothermic surfaces.} Analogously, we arrive at a
(apparently new) characterization of smooth L-isothermic surfaces.
\begin{theorem}\label{thm: L-Darboux}
{\bf (L-isothermic surfaces as Moutard nets in the Laguerre
quadric)} A surface enveloping a two-parameter family of planes
$P:\bbR^2\to\{{\rm planes\ in\ }\bbR^3\}$, with $P=\{x\in\bbR^3:
\langle v,x\rangle=d\}$, is L-isothermic, if and only if there
exists a function $\rho:\bbR^2\to\bbR$ such that
\[
\hat{p}=\rho^{-1}(v+2d\ee_\infty+|v|\ee_6):\bbR^2\to\bbL^{4,2}
\]
satisfies a Moutard equation (\ref{eq: Mou}).
\end{theorem}


\begin{thebibliography}{WWWW}
\bibitem[Do3]{Do3} A.~Doliwa.
   The B-quadrilateral lattice, its transformations and the
   algebro-geometric construction. {\tt arxiv.org/nlin.SI/0601057}.
\end{thebibliography}
\end{document}